\documentclass[12pt]{amsart}

\usepackage[T1]{fontenc}
\usepackage{lmodern}
\usepackage{amsmath,amssymb,amsthm,mathtools}
\usepackage[a4paper,margin=1.05in]{geometry}
\usepackage{microtype}
\usepackage{needspace}
\usepackage{cite}
\usepackage{xcolor}
\usepackage{hyperref}
\hypersetup{colorlinks=true,linkcolor=blue,citecolor=blue,urlcolor=blue,
  pdftitle={Bergman Kernel Zeros on Generalized Hartogs Triangles},
  pdfauthor={Qian Fu}}
\newtheorem{theorem}{Theorem}[section]
\newtheorem*{emquestion}{Edholm--Mathew Question (Q1)}
\newtheorem*{emquestiontwo}{Edholm--Mathew Question (Q2)}

\newtheorem{example}[theorem]{Example}
\newtheorem{lemma}[theorem]{Lemma}
\newtheorem{corollary}[theorem]{Corollary}
\theoremstyle{definition}

\numberwithin{equation}{section}

\makeatletter
\renewenvironment{proof}[1][\proofname]{\par
  \pushQED{\qed}%
  \normalfont \topsep6\p@\@plus6\p@\relax
  \trivlist
  \item[\hskip\labelsep\bfseries #1\@addpunct{.}]\ignorespaces
}{%
  \popQED\endtrivlist\@endpefalse
}
\makeatother

\newcommand{\C}{\mathbb C}
\newcommand{\D}{\mathbb D}

\newcommand{\ceil}[1]{\left\lceil #1\right\rceil}
\newcommand{\floor}[1]{\left\lfloor #1\right\rfloor}

\title[Bergman kernel zeros on generalized Hartogs triangles]
{Bergman Kernel Zeros on Generalized Hartogs Triangles}

\author[Q. Fu]{Qian Fu}
\address{Teaching Experiment Platform, Beijing Normal University, Zhuhai 519087, P. R. China}
\email{qianfu@bnu.edu.cn}

\subjclass[2020]{32A25, 30C15, 42A38}
\keywords{Bergman kernel, generalized Hartogs triangle, Lu Qi-Keng problem}
\date{}

\begin{document}

\begin{abstract}
We study zeros of the Bergman kernel on the generalized Hartogs triangles
$H_\gamma=\{(z_1,z_2)\in\mathbb C^2:|z_1|^\gamma<|z_2|<1\}$.
For every real $1<\gamma<2$, suitable sections near the singular vertex
contain exactly $\lceil(\gamma-1)^{-1}\rceil-1$ zeros in a prescribed
disk, all of them simple. This answers a question of Edholm and Mathew
[\emph{Math. Z.} \textbf{313} (2026)] and, together with the known cases,
shows that $H_\gamma$ has a zero-free
Bergman kernel precisely when $\gamma$ is the reciprocal of a positive
integer. We also fix the second kernel argument and a complex line in
the first argument, and prove that the total number of zeros on the
admissible section is comparable to $(\gamma-1)^{-1}$ as
$\gamma\to1^+$; a comparable number are simple and nonreal.
For the sections near the vertex, all but at most one of the counted
zeros concentrate near the unit circle, and the angular counting error
is uniformly bounded.
\end{abstract}
\maketitle

\section{Introduction}\label{sec:intro}

The Bergman kernel is both the reproducing kernel of a Hilbert space of
holomorphic functions and the integral kernel of an orthogonal
projection. Its positivity on the diagonal does not preclude zeros away
from the diagonal. A bounded domain $\Omega\subset\C^n$ is called a
\emph{Lu Qi-Keng domain} if the reproducing kernel of
$A^2(\Omega)=\mathcal O(\Omega)\cap L^2(\Omega,dV)$ is zero-free on
$\Omega\times\Omega$, where $dV$ denotes Lebesgue volume. This
biholomorphically invariant condition arises in the construction of
Bergman representative coordinates; see Boas~\cite{Boas}.
Explicit examples show that even convex domains can have kernel zeros
\cite{BoasFuStraube1999}. For Hartogs-type domains, explicit formulas
have also connected the zero problem with the roots of auxiliary
polynomials; a detailed treatment for Cartan--Hartogs domains was given
by Ahn and Park~\cite{AhnPark2012}.

We consider the family
\begin{equation}\label{eq:domain}
 H_\gamma=\{(z_1,z_2)\in\C^2:|z_1|^\gamma<|z_2|<1\},
 \qquad \gamma>0,
\end{equation}
and write $K_\gamma=K_{H_\gamma}$. These bounded pseudoconvex Reinhardt
domains provide a model for the effect of boundary singularities on
Bergman theory. Their operator theory is sensitive to the arithmetic of
$\gamma$: the Bergman projection has a nontrivial interval of
$L^p$-boundedness for rational exponents, whereas for irrational
exponents it is bounded only on $L^2$~\cite{EdholmMcNeal2017}.
On the classical triangle, sharp $L^p$ mapping results were obtained by
Chakrabarti and Zeytuncu~\cite{ChakrabartiZeytuncu2016}; weighted norm
estimates and endpoint weak-type behavior were studied by Huo and
Wick~\cite{HuoWickWeighted2020,HuoWickWeak2020}. Related questions of
duality, approximation, and projections onto $L^p$-Bergman spaces have
been developed in \cite{ChakrabartiEdholmMcNeal2019,ChakrabartiEdholm2024}.
The monomial expansions underlying these results will be used here to
study the zeros of the kernel rather than its mapping properties.

For rational exponents, the triangles also admit finite monomial
coverings by $\D\times\D^*$, where
$\D=\{\zeta\in\C:|\zeta|<1\}$ and $\D^*=\D\setminus\{0\}$.
This quotient structure leads to explicit kernel formulas and
transformation laws; see
\cite{Bell1982,BenderChakrabartiEdholmMainkar2022,ChakrabartiKonkelMainkarMiller2020}.
In particular, $(z_1,z_2)\mapsto(z_1/z_2,z_2)$ is a biholomorphism
from $H_1$ onto $\D\times\D^*$, so $K_1$ is zero-free.
Edholm proved that $K_j$ has zeros for every integer $j\ge2$
\cite[Corollary~4.2]{Edholm2016}, and Edholm and McNeal established the
same conclusion for every real $\gamma>2$
\cite[Theorem~3.8]{EdholmMcNeal2017}.
With
\[
 s=z_1\overline{w_1},\qquad t=z_2\overline{w_2},
\]
we also write $K_\gamma(s,t)$ for the holomorphic factor of the kernel.
Its restriction to $s=0$ is
\[
 K_\gamma(0,t)=\frac{1+(\gamma-1)t}{\gamma\pi^2t(1-t)^2},
 \qquad \gamma>1;
\]
see \cite[proof of Theorem~3.8]{EdholmMcNeal2017} and
\cite[Section~1.3.1]{EdholmMathew2026}.
The numerator has a zero in $0<|t|<1$ precisely when $\gamma>2$.
Thus this restriction does not settle the interval $(1,2)$.
Edholm and Mathew~\cite[Section~1]{EdholmMathew2026} formulated the
following question.

\begin{emquestion}
For which $\gamma\in(1,2)$ does the Bergman kernel have zeros?
\end{emquestion}

They also asked whether the zero-free property of the classical
Hartogs triangle persists in a one-sided neighborhood of its exponent.

\begin{emquestiontwo}
Is there some $\varepsilon>0$ so that $K_\gamma(\cdot,\cdot)$ is
non-vanishing on $H_\gamma\times H_\gamma$ for
$\gamma\in[1,1+\varepsilon)$?
\end{emquestiontwo}

Edholm and Mathew answered (Q2) negatively by producing rational
exponents $\gamma_j\to1^+$ for which $K_{\gamma_j}$ has zeros
\cite[Theorem~1.7 and Corollary~1.8]{EdholmMathew2026}.
Their analysis restricts the rational kernel to $s=t$ and introduces a
palindromic polynomial $Q_{m,n}$ when $\gamma=m/n>1$ is in lowest
terms. Edholm and Mathew~\cite[Conjecture~4.20]{EdholmMathew2026}
conjectured that, for every pair of relatively prime positive integers
$m>n$, the associated palindromic polynomial $Q_{m,n}$ has no zeros
on the unit circle and consequently has exactly $m-n$ zeros in the
punctured unit disc, counted with multiplicity. Recently, Deng and
Fu~\cite{DengFu2026} established this conjecture using a
Fej\'er-kernel representation and finite Fourier orthogonality.

In the present paper we instead fix $t$ and vary $s$. The resulting
sections can be studied through a Laurent-series formula valid for
every real exponent. We answer (Q1) throughout $(1,2)$ and obtain an
exact count of simple zeros on sections approaching the singular
vertex. We then quantify the phenomenon behind (Q2) in two ways:
by proving two-sided bounds for the number of zeros on a fixed complex
line, and by describing the angular distribution of the zeros on
sections near the vertex. These are different questions from
unit-circle nonvanishing of the rational polynomials $Q_{m,n}$.

For $1<\gamma<2$ put
\begin{equation}\label{eq:parameters}
 \delta=1-\frac1\gamma,\qquad
 N=\left\lceil\frac1\delta\right\rceil-1,\qquad
 d=N-1=\left\lceil\frac1{\gamma-1}\right\rceil-1 .
\end{equation}
Here $\lceil x\rceil$ denotes the least integer greater than or equal to $x$,
and $\lfloor x\rfloor$ denotes the greatest integer less than or equal to $x$.
For real $0<t<1$ write $s=tv$. The product variables are admissible when $|v|<t^{-\delta}$.
Throughout the paper, a simple zero means a simple zero of the
one-variable function specified in the statement.

\begin{theorem}\label{thm:main}
For each real $1<\gamma<2$ there exist $R_\gamma>3$ and
$t_\gamma\in(0,1)$ such that, for every $0<t<t_\gamma$, the disk
$|v|\le R_\gamma$ lies in $|v|<t^{-\delta}$ and
\[
 v\longmapsto K_\gamma(tv,t)
\]
has exactly $d$ zeros in $|v|<R_\gamma$. All $d$ zeros are simple.
\end{theorem}

The integer $d$ is at least one throughout $(1,2)$. For example, on the
sequence $\gamma=(n+1)/n$, $n\ge2$, Theorem~\ref{thm:main} gives
$n-1$ simple zeros on each sufficiently small positive $t$ section.
For the same exponents, the polynomial associated with the different
restriction $s=t$ has one zero in $\D^*$ by \cite{DengFu2026}.

Together with the known results for $\gamma\ge2$ and the
biholomorphic equivalences of generalized Hartogs triangles,
Theorem~\ref{thm:main} gives a classification for all positive exponents.

\begin{corollary}\label{cor:class}
For $\gamma>0$, $H_\gamma$ is a Lu Qi-Keng domain if and only if
$\gamma=1/j$ for a positive integer $j$.
\end{corollary}

The next theorem keeps both the second kernel argument and a complex
line in the first argument fixed as $\gamma\to1^+$.
It gives matching upper and lower bounds for the total number of zeros
on the admissible part of that line, and shows that a comparable number
are simple and nonreal.

\Needspace{12\baselineskip}
\begin{theorem}\label{thm:fixed}
For each $0<t<3-2\sqrt2$ there exist $a,b\in(0,1)$ with $ab=t$ and
constants $c_t,C_t,\eta_t>0$ such that, for $1<\gamma<1+\eta_t$, the function
\[
 v\longmapsto K_\gamma\bigl((av,a),(b,b)\bigr)
\]
has at least $c_t/(\gamma-1)$ distinct simple nonreal zeros in
$1<|v|<a^{-\delta}$. Its total number $n_\gamma$ of zeros in
$|v|<a^{-\delta}$, counted with multiplicity, satisfies
\begin{equation}\label{eq:fixed-count-intro}
 \frac{c_t}{\gamma-1}\le n_\gamma\le\frac{C_t}{\gamma-1}.
\end{equation}
\end{theorem}

\begin{corollary}\label{cor:fixed-pole}
In Theorem~\ref{thm:fixed} one can take $a=b=1/\sqrt{10}$ for
$t=1/10$. In particular, with $w_*=(1/\sqrt{10},1/\sqrt{10})$, the
number of simple zeros of $z_1\mapsto K_\gamma((z_1,1/\sqrt{10}),w_*)$
in the admissible disk tends to infinity as $\gamma\to1^+$.
\end{corollary}

We return to the sections of Theorem~\ref{thm:main} to study the
location of their zeros as $\gamma\to1^+$.
Angular counting estimates for polynomial roots go back to
Erd\H{o}s and Tur\'an~\cite{ErdosTuran1950}. For sparse polynomials,
Hayman's estimate~\cite{Hayman}, in the form recalled in
\cite[Theorem~1.3]{Hrubes}, bounds the discrepancy in terms of the
number of nonzero coefficients. A sparse representation of the limiting
polynomial will yield a counting error independent of $\gamma$.

\Needspace{7\baselineskip}
Let $I$ be a half-open arc of the unit circle, with length $|I|$
measured in radians; the full circle is also allowed.
For $0<\varepsilon<1$, denote by $Z_{\gamma,t}(\varepsilon,I)$
the number of distinct zeros $v$ of $v\mapsto K_\gamma(tv,t)$
satisfying
\[
 1-\varepsilon<|v|<1+\varepsilon,
 \qquad \frac{v}{|v|}\in I.
\]

\begin{theorem}\label{thm:distribution}
The numbers $t_\gamma$ in Theorem~\ref{thm:main} can be chosen so that
the following holds. For every $0<\varepsilon<1$ there is an
$\eta_\varepsilon>0$ such that, whenever $1<\gamma<1+\eta_\varepsilon$
and $0<t<t_\gamma$, at least $d-1$ of the $d$ zeros lie in
$1-\varepsilon<|v|<1+\varepsilon$. For every half-open arc $I$,
\begin{equation}\label{eq:sector-main}
 \left|Z_{\gamma,t}(\varepsilon,I)-
       \frac{|I|}{2\pi(\gamma-1)}\right|\le 10 .
\end{equation}
\end{theorem}

\section{Preliminary Results}\label{sec:kernel}

We use the usual four-dimensional Lebesgue volume $dV$ on $\mathbb C^2$. For $(z,w)\in H_\gamma\times H_\gamma$, put
\[
 s=z_1\overline{w_1},\qquad t=z_2\overline{w_2}.
\]
Then $0<|t|<1$ and $|s|<|t|^{1/\gamma}$. Conversely, each $(s,t)$ satisfying these two strict inequalities is realized by the choice
\begin{equation}\label{eq:realization}
 z_2=\sqrt{|t|}\,e^{i\arg t},\quad w_2=\sqrt{|t|},\qquad
 z_1=\sqrt{|s|}\,e^{i\arg s},\quad w_1=\sqrt{|s|},
\end{equation}
with $z_1=w_1=0$ when $s=0$. Then $|z_1|^\gamma=|s|^{\gamma/2}<|t|^{1/2}=|z_2|$, and the same holds for $w$.

\begin{lemma}[{\cite[Lemma~3.2]{EdholmMcNeal2017}}]\label{lem:basis}
For $\gamma>0$, the monomials $z_1^a z_2^b$, $a\in\mathbb Z_{\ge0}$, $b\in\mathbb Z$, with
\begin{equation}\label{eq:admissible}
 b+1+\frac{a+1}{\gamma}>0
\end{equation}
form a complete orthogonal basis of $A^2(H_\gamma)$. For each admissible pair,
\begin{equation}\label{eq:norm}
 \left\|z_1^a z_2^b\right\|_{L^2(H_\gamma)}^2
 =\frac{\pi^2}{(a+1)\bigl(b+1+(a+1)/\gamma\bigr)}.
\end{equation}
\end{lemma}

For each positive integer $n$, set
\begin{equation}\label{eq:b-beta}
 \beta_n=1+\frac{n}{\gamma}-\left\lceil\frac{n}{\gamma}\right\rceil
 \in(0,1].
\end{equation}
The integrability condition \eqref{eq:admissible}, with $a=n-1$,
then allows exactly the powers
$z_2^{-\lceil n/\gamma\rceil+k}$, where $k\ge0$ is an integer.
If $n/\gamma$ is an integer, then $\beta_n=1$.

\begin{lemma}\label{lem:kernel}
For every $\gamma>0$ and all product variables $0<|t|<1$, $|s|<|t|^{1/\gamma}$, the Bergman kernel is
\begin{equation}\label{eq:kernel-series}
 K_\gamma(s,t)
 =\frac{1}{\pi^2(1-t)^2}
  \sum_{n=1}^\infty n\,s^{n-1}t^{-\ceil{n/\gamma}}
  \bigl[\beta_n+(1-\beta_n)t\bigr].
\end{equation}
Both the displayed series and the orthogonal-basis expansion used below converge absolutely and locally uniformly on this region.
\end{lemma}

\begin{proof}
For $a=n-1$, the integrability condition in Lemma~\ref{lem:basis}
reads $b>-1-n/\gamma$. Its integer solutions are
\[
 b=-\ceil{n/\gamma}+k,\qquad k=0,1,2,\ldots,
 \qquad b+1+\frac n\gamma=k+\beta_n.
\]
By the reproducing-kernel expansion in the orthogonal basis from
Lemma~\ref{lem:basis}, for $z,w\in H_\gamma$ we have
\begin{align*}
 K_\gamma(z,w)
 &=\sum_{\substack{a\in\mathbb Z_{\ge0},\,b\in\mathbb Z\\
       b+1+(a+1)/\gamma>0}}
   \frac{z_1^az_2^b\overline{w_1^aw_2^b}}
        {\|z_1^az_2^b\|_{L^2(H_\gamma)}^2}\\
 &=\frac1{\pi^2}\sum_{a=0}^\infty
   \sum_{\substack{b\in\mathbb Z\\b+1+(a+1)/\gamma>0}}
   (a+1)\left(b+1+\frac{a+1}{\gamma}\right)
   (z_1\overline{w_1})^a(z_2\overline{w_2})^b.
\end{align*}
We now re-index the sum by setting
\[
 a=n-1,\qquad b=-\ceil{n/\gamma}+k,
 \qquad n\ge1,\quad k\ge0.
\]
The coefficient in the second line becomes
\[
 (a+1)\left(b+1+\frac{a+1}{\gamma}\right)
 =n\left(k+1+\frac n\gamma-\ceil{n/\gamma}\right)
 =n(k+\beta_n).
\]
With $s=z_1\overline{w_1}$ and $t=z_2\overline{w_2}$, this gives
\begin{equation}\label{eq:kernel-double-explicit}
 K_\gamma(s,t)=\frac1{\pi^2}
   \sum_{n=1}^\infty\sum_{k=0}^\infty
   n(k+\beta_n)s^{n-1}t^{-\ceil{n/\gamma}+k}.
\end{equation}
To verify absolute convergence, put
$q=|s|/|t|^{1/\gamma}<1$. Since
$\ceil{n/\gamma}\le n/\gamma+1$ and $0<\beta_n\le1$,
\begin{align*}
\sum_{k=0}^\infty(k+\beta_n)|t|^k
 &=\frac{\beta_n}{1-|t|}+\frac{|t|}{(1-|t|)^2}\\
 &=\frac{\beta_n+(1-\beta_n)|t|}{(1-|t|)^2}
 \le\frac1{(1-|t|)^2}.
\end{align*}
Consequently the absolute-value series in
\eqref{eq:kernel-double-explicit} is bounded by
\begin{align*}
 &\frac1{\pi^2}\sum_{n=1}^\infty
   n|s|^{n-1}|t|^{-\ceil{n/\gamma}}
   \sum_{k=0}^\infty(k+\beta_n)|t|^k\\
 &\qquad\le
 \frac{|t|^{-1-1/\gamma}}{\pi^2(1-|t|)^2}
     \sum_{n=1}^\infty nq^{n-1}
 =\frac{|t|^{-1-1/\gamma}}
       {\pi^2(1-|t|)^2(1-q)^2}.
\end{align*}
On each compact subset of the product-variable region, $|t|$ is bounded away from $0$ and $1$, and $q$ is bounded away from $1$. Thus the bound is locally uniform. We may therefore sum
first in $k$, using
\[
 \sum_{k=0}^\infty(k+\beta_n)t^k
 =\frac{\beta_n}{1-t}+\frac{t}{(1-t)^2}
 =\frac{\beta_n+(1-\beta_n)t}{(1-t)^2}.
\]
Substitution into \eqref{eq:kernel-double-explicit} proves the formula and its local uniform convergence.
\end{proof}

For instance, setting $s=0$ in \eqref{eq:kernel-series} gives
\begin{equation}\label{eq:s0}
 K_\gamma(0,t)=\frac{1+(\gamma-1)t}{\gamma\pi^2t(1-t)^2}
 \qquad (\gamma>1),
\end{equation}
as in \cite[proof of Theorem~3.8]{EdholmMcNeal2017}.

We next normalize the kernel near the singular vertex. For $1<\gamma<2$, recall that
\begin{equation}\label{eq:parameters-recall}
 \delta=1-\gamma^{-1}\in(0,\tfrac12),
 \qquad N=\ceil{1/\delta}-1\ge2.
\end{equation}
Let $0<t<1$ be real. For $|v|<t^{-\delta}$, the product variables $(s,t)=(vt,t)$ satisfy $|vt|^\gamma<t$. Define
\begin{equation}\label{eq:Fdef}
 F_t(v)=\pi^2t(1-t)^2K_\gamma(vt,t).
\end{equation}
Here $\{x\}=x-\floor{x}$ denotes the fractional part of $x$.
Since $n/\gamma=n-\delta n$, we have
\[
 n-\ceil{n/\gamma}=\floor{\delta n},\qquad
 \beta_n=1-\{\delta n\}.
\]
Substituting these identities into \eqref{eq:kernel-series} gives
\begin{equation}\label{eq:section}
 F_t(v)=\sum_{n=1}^\infty n\,v^{n-1}
   t^{\floor{\delta n}}
   \bigl[\beta_n+(1-\beta_n)t\bigr].
\end{equation}
All coefficients of this series are real when $t>0$.

For fixed $\gamma$, the factor $t^{\floor{\delta n}}$ in
\eqref{eq:section} singles out the terms that survive as
$t\to0^+$.  The following lemma shows that only finitely many
terms survive and gives an explicit bound for the remaining series.

\begin{lemma}\label{lem:degeneration}
Fix $1<\gamma<2$ and let $\delta,N$ be as in
\eqref{eq:parameters}.  Then the following statements hold.
\begin{enumerate}
\renewcommand{\labelenumi}{\textup{(\roman{enumi})}}
\renewcommand{\theenumi}{(\roman{enumi})}
\item\label{lem:degeneration:i} The integer $N$ is the largest positive integer satisfying
$\delta N<1$.  The polynomial
\begin{equation}\label{eq:poly}
 P_\delta(v)=\sum_{n=1}^{N}n(1-n\delta)v^{n-1}
\end{equation}
has degree $N-1\ge1$ and strictly positive coefficients.

\item\label{lem:degeneration:ii} As $t\to0^+$, the domains $|v|<t^{-\delta}$ expand to
$\mathbb C$, and the normalized kernel sections satisfy
\[
 F_t(v)\longrightarrow P_\delta(v)
\]
uniformly on every compact subset of $\mathbb C$.

\item\label{lem:degeneration:iii} More precisely, fix $R>0$ and put
\begin{equation}\label{eq:hq}
 q=Rt^{\delta/2}.
\end{equation}
If $0<t<1$ and $q<1$, then $|v|\le R$ lies in the domain of
$F_t$, and
\begin{equation}\label{eq:error}
 \sup_{|v|\le R}|F_t(v)-P_\delta(v)|
 \le t\delta\sum_{n=1}^{N}n^2R^{n-1}
 +t^{\delta/2}\frac{q^{N}\bigl((N+1)-Nq\bigr)}{(1-q)^2}.
\end{equation}
\end{enumerate}
\end{lemma}

\Needspace{9\baselineskip}
\begin{proof}
The definition of the ceiling function and
$N=\ceil{1/\delta}-1$ give
\[
 N<\frac1\delta\le N+1.
\]
Since $\delta>0$, multiplication by $\delta$ yields
\begin{equation}\label{eq:cutoff}
 \delta N<1\le\delta(N+1).
\end{equation}
For every positive integer $n$ there are two possibilities:
\[
\begin{aligned}
 1\le n\le N
 &\ \Longrightarrow\ 
 0<\delta n\le\delta N<1
 \ \Longrightarrow\ \floor{\delta n}=0,\\
 n\ge N+1
 &\ \Longrightarrow\
 \delta n\ge\delta(N+1)\ge1
 \ \Longrightarrow\ \floor{\delta n}\ge1.
\end{aligned}
\]
Consequently,
\begin{equation}\label{eq:cutoff-indices}
 \floor{\delta n}=0
 \quad\Longleftrightarrow\quad
 1\le n\le N,\qquad n\in\mathbb Z_{>0}.
\end{equation}
For these first $N$ indices, $\{\delta n\}=\delta n$ and hence
$\beta_n=1-\delta n$.  Their contribution to
\eqref{eq:section} is therefore
\begin{align*}
 &\sum_{n=1}^{N}n v^{n-1}
       t^{\floor{\delta n}}[\beta_n+(1-\beta_n)t]\\
 &\quad=\sum_{n=1}^{N}n(1-n\delta)v^{n-1}
         +t\delta\sum_{n=1}^{N}n^2v^{n-1}\\
 &\quad=P_\delta(v)+t\delta\sum_{n=1}^{N}n^2v^{n-1}.
\end{align*}
By \eqref{eq:cutoff}, $1-\delta n>0$ for all
$1\le n\le N$.  Thus every coefficient of $P_\delta$ is
strictly positive, and its leading coefficient
$N(1-N\delta)$ is nonzero.  It follows that
$\deg P_\delta=N-1$. This proves part~\ref{lem:degeneration:i}.

It remains to control the terms indexed by $n\ge N+1$.
Set $j=\floor{\delta n}$; by \eqref{eq:cutoff-indices},
$j\ge1$.  Since $j\le\delta n<j+1\le2j$,
\[
 \floor{\delta n}=j\ge\frac{\delta n}{2},
 \qquad
 0<t^{\floor{\delta n}}\le t^{\delta n/2}
 \quad (0<t<1).
\]
Furthermore, $0<\beta_n\le1$ implies
$0<\beta_n+(1-\beta_n)t\le1$.  Therefore
\begin{align*}
 &\sum_{n=N+1}^\infty nR^{n-1}
       t^{\floor{\delta n}}[\beta_n+(1-\beta_n)t]\\
 &\quad\le\sum_{n=N+1}^\infty nR^{n-1}t^{\delta n/2}\\
 &\quad=t^{\delta/2}\sum_{n=N+1}^\infty nq^{n-1}
 =t^{\delta/2}\frac{q^N((N+1)-Nq)}{(1-q)^2},
 \qquad q=Rt^{\delta/2}.
\end{align*}
Combining the finite part and the tail proves \eqref{eq:error}.
Moreover,
\[
 Rt^\delta=(Rt^{\delta/2})t^{\delta/2}
    =qt^{\delta/2}<1,
\]
so the full disk $|v|\le R$ lies in the section's domain.
This proves part~\ref{lem:degeneration:iii}.

Fix $R>0$. Since $q=Rt^{\delta/2}\to0$ as $t\to0^+$,
we may assume that $0<q\le\frac12$. It follows that
\[
 0<\frac{(N+1)-Nq}{(1-q)^2}\le4(N+1).
\]
Using this bound in \eqref{eq:error}, we obtain
\begin{align*}
 0\le\sup_{|v|\le R}|F_t(v)-P_\delta(v)|
 &\le t\delta\sum_{n=1}^{N}n^2R^{n-1}
    +4(N+1)t^{\delta/2}q^N\\
 &=t\delta\sum_{n=1}^{N}n^2R^{n-1}
    +4(N+1)R^N t^{\delta(N+1)/2}.
\end{align*}
Here $\delta$, $N$, and $R$ are fixed, and
$\delta(N+1)\ge1$ by \eqref{eq:cutoff}. Consequently,
\[
\begin{aligned}
 \lim_{t\to0^+}t\delta\sum_{n=1}^{N}n^2R^{n-1}&=0,\\
 \lim_{t\to0^+}4(N+1)R^N t^{\delta(N+1)/2}&=0.
\end{aligned}
\]
Therefore,
\[
 \lim_{t\to0^+}\sup_{|v|\le R}|F_t(v)-P_\delta(v)|=0.
\]
Finally, $t^{-\delta}\to\infty$, so every compact subset of
$\mathbb C$ lies in the domain of $F_t$ for all sufficiently small
$t>0$. This proves part~\ref{lem:degeneration:ii}.
\end{proof}

\begin{lemma}[Rouch\'e's theorem {\normalfont\cite[p.~153]{Ahlfors1979}}]\label{lem:rouche}
Let $D=\{z\in\C:|z-z_0|<r\}$, where $r>0$, and let $f$ and
$g$ be holomorphic on a neighborhood of $\overline D$. If
\[
 |f(z)-g(z)|<|g(z)| \qquad (|z-z_0|=r),
\]
then $f$ and $g$ have the same number of zeros in $D$, counted
with multiplicity.
\end{lemma}

\section{Proof of Theorem~\ref{thm:main} and Corollary~\ref{cor:class}}
\label{sec:main}

The polynomial limit gives the number of zeros with multiplicity.
To prove simplicity for every sufficiently small positive $t$, we first identify the possible multiple roots of the limit and then compute the first correction to the kernel section at such a root.

\begin{lemma}\label{lem:multiplicity}
Let $1<\gamma<2$, put $\delta=1-\gamma^{-1}$, and define
$N=\ceil{\delta^{-1}}-1\ge2$ and $d=N-1$.
The polynomial $P_\delta$ in Lemma~\ref{lem:degeneration}\ref{lem:degeneration:i}
has at least $d-1$ distinct zeros. All its zeros are simple,
with the possible exception of one double zero on the negative
real axis.
\end{lemma}

\begin{proof}
By the definition of the ceiling function,
\[
 N=\ceil{\delta^{-1}}-1
 \quad\Longrightarrow\quad
 N<\delta^{-1}\le N+1.
\]
Consequently, we can write
\[
 \delta^{-1}=N+\alpha,
 \qquad \alpha:=\delta^{-1}-N,
 \qquad 0<\alpha\le1.
\]
Since $1-n\delta=\delta(N+\alpha-n)$, the polynomial
in Lemma~\ref{lem:degeneration}\ref{lem:degeneration:i} satisfies
\[
 \frac{P_\delta(z)}{\delta}
 =\sum_{n=1}^{N}n(N+\alpha-n)z^{n-1}.
\]
For $|z|<1$, the geometric-series identities give
\[
 \sum_{n=1}^{\infty}nz^{n-1}=\frac1{(1-z)^2},
 \qquad
 \sum_{n=1}^{\infty}n^2z^{n-1}=\frac{1+z}{(1-z)^3}.
\]
Thus
\begin{align*}
 \sum_{n=1}^{\infty}n(N+\alpha-n)z^{n-1}
 &=(N+\alpha)\sum_{n=1}^{\infty}nz^{n-1}
      -\sum_{n=1}^{\infty}n^2z^{n-1}\\
 &=\frac{N+\alpha}{(1-z)^2}
      -\frac{1+z}{(1-z)^3}\\
 &=\frac{N+\alpha-1-(N+\alpha+1)z}{(1-z)^3}.
\end{align*}
To recover the finite sum, we subtract the terms with $n\ge N+1$.
Putting $n=N+1+\ell$ in the tail gives
\begin{align*}
 &\sum_{n=N+1}^{\infty}n(N+\alpha-n)z^{n-1}\\
 &\quad=z^N\sum_{\ell=0}^{\infty}
             (N+1+\ell)(\alpha-1-\ell)z^\ell\\
 &\quad=z^N\left[
       (N+1)(\alpha-1)\sum_{\ell=0}^{\infty}z^\ell
      +(\alpha-N-2)\sum_{\ell=0}^{\infty}\ell z^\ell
      -\sum_{\ell=0}^{\infty}\ell^2z^\ell\right]\\
 &\quad=\frac{z^N}{(1-z)^3}
  \left[(N+1)(\alpha-1)(1-z)^2
       +(\alpha-N-2)z(1-z)-z(1+z)\right]\\
 &\quad=-\frac{z^N}{(1-z)^3}
 \left[(N+1)(1-\alpha)
       +(2N\alpha+\alpha-N+1)z-N\alpha z^2\right].
\end{align*}
Subtracting this tail from the infinite series and multiplying
by $(1-z)^3$ gives
\begin{align}\label{eq:sparse-limit}
 S(z)&:=\frac{(1-z)^3}{\delta}P_\delta(z)
      =A(z)+z^N B(z),\\
 A(z)&=N+\alpha-1-(N+\alpha+1)z,\notag\\
 B(z)&=(N+1)(1-\alpha)
     +(2N\alpha+\alpha-N+1)z-N\alpha z^2.\notag
\end{align}
The calculation was carried out for $|z|<1$, but both sides of
\eqref{eq:sparse-limit} are polynomials; hence the identity holds
for every $z\in\C$.

Differentiating $S=A+z^NB$, we obtain
\[
 S'(z)=A'(z)+Nz^{N-1}B(z)+z^N B'(z).
\]
The terms containing $z^N$ cancel in the combination
\begin{align*}
 &zB(z)S'(z)-\bigl(NB(z)+zB'(z)\bigr)S(z)\\
 &\quad=\bigl(zA'(z)-NA(z)\bigr)B(z)-zA(z)B'(z)\\
 &\quad=-N(N+1)(z-1)^2
 \bigl[\alpha(N+\alpha+1)z
       +(1-\alpha)(N+\alpha-1)\bigr].
\end{align*}
Suppose that $w$ is a zero of $P_\delta$ of multiplicity $r\ge2$.
Since
\[
 P_\delta(0)=1-\delta>0,
 \qquad P_\delta(1)=\sum_{n=1}^Nn(1-n\delta)>0,
\]
we have $w\notin\{0,1\}$.  It follows that $S(w)=S'(w)=0$.
The preceding identity therefore implies
\[
 \alpha(N+\alpha+1)w
       +(1-\alpha)(N+\alpha-1)=0,
 \qquad
 w=-\frac{(1-\alpha)(N+\alpha-1)}
          {\alpha(N+\alpha+1)}.
\]
If $\alpha=1$, then $w=0$, a contradiction. Thus
$0<\alpha<1$ and $w<0$. The displayed formula also shows that
there can be at most one multiple zero.

Finally, if $r\ge3$, then $S$ and $S'$ vanish at $w$ to
orders at least $r$ and $r-1\ge2$, respectively. The left-hand
side of the elimination identity would have a zero of order at
least two, while the right-hand side has a simple zero at
$w\ne1$, because its last factor is linear with nonzero slope
$\alpha(N+\alpha+1)$. This contradiction proves that $r=2$.
\end{proof}

The possible double zero already occurs in degree two. When $N=3$,
\[
 P_\delta(z)/\delta=(2+\alpha)+2(1+\alpha)z+3\alpha z^2 .
\]
Its discriminant is $4(1-4\alpha-2\alpha^2)$, which vanishes at
$\alpha=(\sqrt6-2)/2$, or equivalently $\gamma=\sqrt6-1$.

\begin{lemma}\label{lem:splitting}
Suppose that $P_\delta$ has a double zero $w$. For all sufficiently
small $t>0$, the section $F_t$ has two simple negative real zeros
converging to $w$. More precisely, they can be labeled $w_-(t),w_+(t)$
so that
\[
 w_\pm(t)=w\pm
 \sqrt{-\frac{4P_{\delta/2}(w)}{P_\delta''(w)}}\,\sqrt t
 +o(\sqrt t).
\]
\end{lemma}

\begin{proof}
We retain the polynomials $A,B$ from Lemma~\ref{lem:multiplicity} and the identity
\eqref{eq:sparse-limit}. By Lemma~\ref{lem:multiplicity},
$w$ is a double root and
\begin{equation}\label{eq:double-w}
 w=-\frac{(1-\alpha)(N+\alpha-1)}{\alpha(N+\alpha+1)},
 \qquad 0<\alpha<1.
\end{equation}
Substituting into $A$ and $B$ yields
\begin{align}
 A(w)&=\frac{N+\alpha-1}{\alpha},\notag\\
 B(w)&=\frac{1-\alpha}{\alpha}
      \left(1-\frac{4N(1-\alpha)}{(N+\alpha+1)^2}\right).
 \label{eq:double-AB}
\end{align}
Since $(N+\alpha+1)^2>(N+1)^2\ge4N$,
\[
 0<B(w)<\frac{1-\alpha}{\alpha}<\frac{N+\alpha-1}{\alpha}=A(w).
\]
The equality $S(w)=0$ now gives
\begin{equation}\label{eq:double-power}
 w^N=-\frac{A(w)}{B(w)}<-1.
\end{equation}
As $w<0$, necessarily $N$ is odd and $w<-1$; in particular
$N\ge3$. Furthermore, \eqref{eq:double-w} implies
\[
 (1-\alpha)(N+\alpha-1)>\alpha(N+\alpha+1)
 \quad\Longrightarrow\quad
 2\alpha(N+\alpha)<N+\alpha-1
 \quad\Longrightarrow\quad 0<\alpha<\tfrac12.
\]
It follows that $2(N+\alpha)=2N+2\alpha$, where $0<2\alpha<1$; the cutoff
index for $P_{\delta/2}$ is therefore $2N$.

Let $A_2,B_2$ be the polynomials in \eqref{eq:sparse-limit}
with $(N,\alpha)$ replaced by $(2N,2\alpha)$. Substituting the expression for $w$ in \eqref{eq:double-w} gives
\begin{align}\label{eq:double-B-identity}
 &-\alpha(N+\alpha+1)^2\bigl(B_2(w)+B(w)\bigr)\notag\\
 &\quad=2\bigl[N(N-2)^2-1\bigr]
       +\alpha(17N^2-36N-3)+(40N+4)\alpha^2+5\alpha^3>0.
\end{align}
All coefficients on the right are positive when $N\ge3$.
Also
\[
 A_2(w)=2A(w)+1+w<2A(w),
 \qquad B_2(w)<-B(w),
\]
and \eqref{eq:double-AB} gives
\[
 \frac{A(w)}{B(w)}>
   \frac{N+\alpha-1}{1-\alpha}>2.
\]
Combining these inequalities with
$w^{2N}=A(w)^2/B(w)^2$, we obtain
\begin{align*}
 \frac{(1-w)^3}{\delta/2}P_{\delta/2}(w)
 &=A_2(w)+w^{2N}B_2(w)\\
 &<2A(w)-\frac{A(w)^2}{B(w)^2}B(w)\\
 &=A(w)\left(2-\frac{A(w)}{B(w)}\right)<0.
\end{align*}
Thus $P_{\delta/2}(w)<0$.

We next determine the sign of $P_\delta''(w)$. Differentiating
the elimination identity in Lemma~\ref{lem:multiplicity} and
using $S(w)=S'(w)=0$, we get
\[
 wB(w)S''(w)
 =-N(N+1)(w-1)^2 \alpha(N+\alpha+1)<0.
\]
Since $wB(w)<0$, this implies $S''(w)>0$; at the double root,
\[
 S''(w)=\frac{(1-w)^3}{\delta}P_\delta''(w),
 \qquad\text{so}\qquad P_\delta''(w)>0.
\]

It remains to compare the normalized kernel section with its limit
one order beyond \eqref{eq:poly}. For $1\le n\le N$, its coefficient
in \eqref{eq:section} is
\[
 n\bigl(1-n\delta+n\delta t\bigr),
\]
whereas for $1\le\delta n<2$ the coefficient of $t$ is
$n(2-n\delta)$. All remaining terms have order at least $t^2$.
For each fixed $R>0$, choose $\rho\in(0,1)$ with
$R\rho^\delta<1$. The series
\[
 \sum_{n\ge1}nR^{n-1}\rho^{\floor{\delta n}}
 \le \rho^{-1}\sum_{n\ge1}nR^{n-1}\rho^{\delta n}
 <\infty
\]
provides a summable majorant on $|v|\le R$, $|t|\le\rho$. Thus the series extends holomorphically to $t=0$, and its Taylor expansion gives
\begin{equation}\label{eq:second-taylor}
 F_t(v)=P_\delta(v)
       +2t\bigl(P_{\delta/2}(v)-P_\delta(v)\bigr)
       +O_R(t^2),\qquad |v|\le R.
\end{equation}
Set $v=w+\sqrt t\,u$. Using $P_\delta(w)=P_\delta'(w)=0$,
\begin{align*}
 P_\delta(w+\sqrt t\,u)
 &=\tfrac12 tP_\delta''(w)u^2+O(t^{3/2}),\\
 2t\bigl(P_{\delta/2}(w+\sqrt t\,u)
       -P_\delta(w+\sqrt t\,u)\bigr)
 &=2tP_{\delta/2}(w)+O(t^{3/2}).
\end{align*}
Hence, locally uniformly in $u$,
\[
 t^{-1}F_t(w+\sqrt t\,u)
 \longrightarrow \tfrac12P_\delta''(w)u^2+2P_{\delta/2}(w).
\]
The limiting quadratic has the two distinct real zeros
\[
 u_\pm=\pm\sqrt{-\frac{4P_{\delta/2}(w)}{P_\delta''(w)}}.
\]
Applying Lemma~\ref{lem:rouche} to disjoint disks centered at $u_-$ and $u_+$
shows that the scaled section has exactly one zero in each disk
for all small $t>0$. Its coefficients are real, so uniqueness in
each conjugation-invariant disk forces both zeros to be real.
Letting the disk radii tend to zero gives
\[
 w_\pm(t)=w+\sqrt t\,u_\pm+o(\sqrt t).
\]
Both are negative for small $t$, because $w<-1$.
\end{proof}

\begin{proof}[Proof of Theorem~\ref{thm:main}]
Recall
\[
 P_\delta(v)=\sum_{n=1}^{N}n(1-n\delta)v^{n-1},
 \qquad F_t(v)=\pi^2t(1-t)^2K_\gamma(tv,t).
\]
Choose
\[
 R_\gamma=4+
   \frac{\displaystyle\sum_{n=1}^{N-1}n(1-n\delta)}{N(1-N\delta)}>3.
\]
Since $R_\gamma>1$, on $|v|=R_\gamma$,
\begin{align*}
 &N(1-N\delta)|v|^{N-1}
  -\left|\sum_{n=1}^{N-1}n(1-n\delta)v^{n-1}\right|\\
 &\quad\ge N(1-N\delta)R_\gamma^{N-1}
   -R_\gamma^{N-2}\sum_{n=1}^{N-1}n(1-n\delta)\\
 &\quad=4N(1-N\delta)R_\gamma^{N-2}>0.
\end{align*}
Lemma~\ref{lem:rouche} shows that $P_\delta$ has all
$d=N-1$ zeros, counted with multiplicity, in $|v|<R_\gamma$.
Since $P_\delta(0)=1-\delta>0$, none is at the origin.

Choose pairwise disjoint closed disks about its distinct roots,
contained in $0<|v|<R_\gamma$, with no polynomial root on a
boundary. On the union of these boundaries and $|v|=R_\gamma$,
the continuous function $|P_\delta|$ has a positive minimum.
By Lemma~\ref{lem:degeneration}\ref{lem:degeneration:ii}, as $t\to0^+$,
\[
 \sup_{|v|\le R_\gamma}|F_t(v)-P_\delta(v)|\longrightarrow0.
\]
Consequently, after choosing $t_\gamma>0$ sufficiently small,
Lemma~\ref{lem:rouche} applies on every boundary. We also impose
$R_\gamma t_\gamma^\delta<1$, so these disks lie in the domain
of the normalized section. Thus $F_t$ has exactly $d$ zeros in
$|v|<R_\gamma$, with the same number in each small disk as the
multiplicity of its center.

For every simple root of $P_\delta$, the corresponding disk
contains one simple root of $F_t$. Lemma~\ref{lem:multiplicity}
allows at most one double root; Lemma~\ref{lem:splitting} shows
that its two nearby roots of $F_t$ are also simple. Decreasing
$t_\gamma$ once more if necessary establishes simplicity for all
$0<t<t_\gamma$.

Finally, on $|v|\le R_\gamma$,
\[
 |s|=t|v|\le tR_\gamma<t^{1/\gamma}
 \quad\Longleftrightarrow\quad R_\gamma t^\delta<1.
\]
The factor $\pi^2t(1-t)^2$ never vanishes. Every zero of $F_t$
therefore gives an actual Bergman kernel zero through
\eqref{eq:realization}.
\end{proof}

\begin{proof}[Proof of Corollary~\ref{cor:class}]
\enlargethispage{2\baselineskip}
For $1<\gamma<2$ the kernel has zeros by Theorem~\ref{thm:main};
for $\gamma\ge2$ this follows from
\cite{Edholm2016,EdholmMcNeal2017}. We now consider $0<\gamma\le1$.

If $\gamma=1/j$ with $j\in\mathbb Z_+$, define
\[
 \Psi_j(z_1,z_2)=(z_1z_2^{-j},z_2),\qquad
 \Psi_j^{-1}(u,v)=(uv^j,v).
\]
The equivalence
\[
 |z_1|<|z_2|^j<1
 \quad\Longleftrightarrow\quad
 |z_1z_2^{-j}|<1,\quad 0<|z_2|<1
\]
shows that $\Psi_j$ is a biholomorphism
$H_{1/j}\simeq\D\times\D^*$ with
$\det\Psi_j'(z)=z_2^{-j}\ne0$.
For $f(v)=\sum_{k\in\mathbb Z}c_kv^k\in A^2(\D^*)$,
angular orthogonality implies, for any $k\le-1$ and
$0<\varepsilon<1/2$,
\[
 \|f\|_2^2\ge2\pi|c_k|^2
   \int_\varepsilon^{1/2}r^{2k+1}\,dr
 \xrightarrow[\varepsilon\to0^+]{}\infty
 \quad\text{unless }c_k=0.
\]
Thus $A^2(\D^*)=A^2(\D)$ isometrically, and
\[
 K_{\D\times\D^*}((u,v),(\xi,\eta))
 =\frac1{\pi^2(1-u\bar\xi)^2(1-v\bar\eta)^2}\ne0.
\]
The Bergman kernel transformation law \cite{Bell1982} therefore shows that $K_{1/j}$ is zero-free.

Otherwise $0<\gamma<1$ and $\gamma^{-1}\notin\mathbb Z$.
Write
\[
 \gamma^{-1}=j+\tau,\qquad
 j=\floor{\gamma^{-1}}\ge1,\quad 0<\tau<1.
\]
The maps
\[
 \Phi(z_1,z_2)=(z_1z_2^{-j},z_2),\qquad
 \Phi^{-1}(u_1,u_2)=(u_1u_2^j,u_2)
\]
define inverse biholomorphisms between $H_\gamma$ and $H_{1/\tau}$ \cite[Theorem~2.19]{EdholmMathew2026}, since
\begin{align*}
 |z_1|^\gamma<|z_2|<1
 &\ \Longleftrightarrow\ |z_1|<|z_2|^{j+\tau}<1\\
 &\ \Longleftrightarrow\ |z_1z_2^{-j}|<|z_2|^\tau<1.
\end{align*}
In particular $\det\Phi'(z)=z_2^{-j}\ne0$.
Since $1/\tau>1$, the Bergman kernel of $H_{1/\tau}$
has zeros by the first part of the proof. For points whose images
are such a zero pair, the transformation law gives
\[
 K_\gamma(z,w)
 =z_2^{-j}K_{1/\tau}(\Phi(z),\Phi(w))\,
       \overline{w_2^{-j}}=0.
\]
This leaves exactly the zero-free exponents $\gamma=1/j$.
\end{proof}

We finish this section with an explicit illustration of Theorem~\ref{thm:main} for the irrational exponent $\gamma=\sqrt2$.

\begin{example}\label{prop:sqrt2-intro}
Set $\gamma=\sqrt2$, $t=2^{-40}$, and
\begin{equation}\label{eq:centers}
 v_\pm=-\frac{2+\sqrt2}{3}\ \pm\ i\frac{1+\sqrt2}{3}.
\end{equation}
The function $v\mapsto K_{\sqrt2}(vt,t)$ has exactly one simple zero in each of the disjoint disks $|v-v_\pm|<\frac14$. Both zeros are nonreal and conjugate. They are the only zeros, counted with multiplicity, in $|v|<4$, and all points of these disks are admissible product variables for $H_{\sqrt2}$.
\end{example}

\begin{proof}
Write $r=\sqrt2$ and $\delta=1-1/r$. Since $1/4<\delta<1/3$, we have $N=3$. The polynomial from Lemma~\ref{lem:degeneration}\ref{lem:degeneration:i} becomes
\begin{equation}\label{eq:quadratic}
 P_\delta(v)=\frac1r+(2r-2)v+cv^2,
 \qquad c=\frac9r-6>0.
\end{equation}
Its discriminant is
\[
 (2r-2)^2-\frac{4c}{r}
 =4r-6=-(2-r)^2,
\]
so its roots are precisely the centers $v_\pm$ in \eqref{eq:centers} and
\begin{equation}\label{eq:quadratic-factor}
 P_\delta(v)=c(v-v_+)(v-v_-).
\end{equation}
The inequalities $27>19r$ and $r>5/4$ give
\begin{equation}\label{eq:sep}
 c>\frac13,\qquad
 |v_+-v_-|=\frac{2(1+r)}3>\frac32.
\end{equation}
Moreover $|v_\pm|^2=1+2r/3<2$, so the two radius-$1/4$ disks are disjoint, lie inside $|v|<7/4<4$, and do not meet the real axis. On each boundary circle $|v-v_\pm|=1/4$, equations \eqref{eq:quadratic-factor}--\eqref{eq:sep} imply
\begin{equation}\label{eq:small-P}
 |P_\delta(v)|\ge\frac c4\left(|v_+-v_-|-\frac14\right)
 >\frac13\cdot\frac14\cdot\frac54=\frac5{48}.
\end{equation}

We estimate $F_t-P_\delta$ on $|v|\le4$ with $t=2^{-40}$. Since $t^{\delta/2}=2^{-20\delta}<1/32$,
we also have $4t^{\delta/2}<1/8$. By Lemma~\ref{lem:degeneration}\ref{lem:degeneration:iii}, the finite correction is bounded by
\[
 t\delta\sum_{n=1}^3n^2 4^{n-1}
 =161\delta\,2^{-40}<161\,2^{-40}<\frac1{3000}.
\]
The infinite remainder is strictly smaller than
\[
 \frac1{32}\cdot \left(\frac18\right)^3
 \frac4{(1-1/8)^2}=\frac1{3136}<\frac1{3000}.
\]
Consequently,
\begin{equation}\label{eq:sqrt-error}
 \sup_{|v|\le4}|F_t(v)-P_\delta(v)|
 <\frac1{1500}<\frac5{48}.
\end{equation}
Moreover,
\[
 4t^\delta<4\cdot2^{-10}=1/256<1,
\]
so $F_t$ is holomorphic near the closed disk.

Applying Lemma~\ref{lem:rouche} on each of the small circles shows that $F_t$ has exactly one zero, counted with multiplicity, in each disk, because the polynomial has exactly one there. In particular the zeros are simple. They are nonreal because both disks avoid the real axis, and they are conjugate because \eqref{eq:section} has real coefficients for real $t$.

To count all zeros in $|v|<4$, observe that on $|v|=4$ the leading term in \eqref{eq:quadratic} dominates the lower terms by
\[
 16c-4(2r-2)-\frac1r=\frac{127}{r}-88>1.
\]
The last inequality is equivalent to $127>89r$ and follows from $127^2=16129>15842=2\cdot89^2$. Combining Lemma~\ref{lem:rouche} with \eqref{eq:sqrt-error} shows that $F_t$ has exactly two zeros in $|v|<4$. The two already localized zeros therefore exhaust them.

Finally, $|v|t^\delta<1$ on the closed disk, so each zero gives an admissible pair $(s,t)=(vt,t)$, and hence actual kernel-zero points of $H_{\sqrt2}$ by \eqref{eq:realization}.
\end{proof}

\section{Proof of Theorem~\ref{thm:fixed}}\label{sec:fixed}

Fix $0<t<3-2\sqrt2$ and set $L=-\log t$.
To keep the dependence on both parameters visible, write
\[
 F_{\delta,t}(v)=\pi^2t(1-t)^2 K_\gamma(tv,t),
 \qquad \delta=1-\gamma^{-1}.
\]
The coefficients in the exact series are obtained by sampling the following
continuous, piecewise linear function:
\[
 h_t(x)=t^{\lfloor x\rfloor}
   \bigl(1-(1-t)(x-\lfloor x\rfloor)\bigr),\qquad x\ge0.
\]
Indeed, the series for the rescaled kernel gives
\begin{equation}\label{eq:fixed-samples}
 F_{\delta,t}(v)=\sum_{n=1}^\infty n h_t(\delta n)v^{n-1},
 \qquad |v|<t^{-\delta}.
\end{equation}
We use the logarithmic variable $v=e^{\delta u}$, for which the
domain of convergence in \eqref{eq:fixed-samples} is $\Re u<L$.

\begin{lemma}\label{lem:fixed-periodization}
Define
\[
 H_t(u)=\int_0^\infty x h_t(x)e^{ux}\,dx,\qquad \Re u<L .
\]
For each compact interval $J\subset(-\infty,L)$ there is a constant
$C_J$, depending also on $t$, such that
\begin{equation}\label{eq:fixed-approx}
 \left|\delta^2 e^{\delta u}F_{\delta,t}(e^{\delta u})-H_t(u)\right|
 \le C_J\delta^2,\qquad
 \Re u\in J,\quad |\Im u|\le\pi/\delta,\quad 0<\delta<1 .
\end{equation}
\end{lemma}

\begin{proof}
For $k\le x<k+1$, $k\ge0$, the sampling function is
\[
 h_t(x)=t^k\bigl(1-(1-t)(x-k)\bigr),\qquad
 h_t(k)=t^k,
\]
so it is continuous at each integer and has slope
$h_t'(x)=-(1-t)t^k$ on $(k,k+1)$. Moreover,
\[
 0<h_t(x)\le t^k\le t^{-1}e^{-Lx},
 \qquad L=-\log t.
\]
Hence all integrals below converge absolutely for $\Re u<L$.
Set
$J_t(u)=\int_0^\infty h_t(x)e^{ux}\,dx$.
For $u\ne0$, integration by parts gives
\begin{align*}
 uJ_t(u)
 &=[h_t(x)e^{ux}]_0^\infty
        -\int_0^\infty h_t'(x)e^{ux}\,dx\\
 &=-1+(1-t)\sum_{k=0}^\infty t^k
        \int_k^{k+1}e^{ux}\,dx\\
 &=-1+\frac{(1-t)(e^u-1)}u
               \sum_{k=0}^\infty(te^u)^k\\
 &=-1+\frac{(1-t)(e^u-1)}{u(1-te^u)}.
\end{align*}
Thus
\[
 J_t(u)=-\frac1u
    +\frac{(1-t)(e^u-1)}{u^2(1-te^u)}.
\]
The apparent singularity at $u=0$ is removable, because the
integral defining $J_t$ is holomorphic there. Differentiation
under the integral sign yields $H_t=J_t'$. In particular,
\begin{equation}\label{eq:fixed-transform}
\begin{aligned}
 u^2H_t(u)
 &=1+\frac{d}{du}\left(\frac{(1-t)(e^u-1)}{1-te^u}\right)
    -\frac{2(1-t)(e^u-1)}{u(1-te^u)}\\
 &=1+\frac{(1-t)^2e^u}{(1-te^u)^2}
    -\frac{2(1-t)(e^u-1)}{u(1-te^u)}.
\end{aligned}
\end{equation}

For $\Re u\in J\Subset(-\infty,L)$,
\[
 |1-te^u|\ge1-te^{\sup J}>0,
 \qquad |e^u|\le e^{\sup J}.
\]
Thus $|u^2H_t(u)|$ is bounded when $|\Im u|\ge1$;
when $|\Im u|\le1$, the integral defining $H_t$ is uniformly
bounded on the corresponding compact strip. Consequently
\[
 |H_t(u)|\le \frac{C_J}{1+|\Im u|^2}.
\]

For the Poisson summation identity, cf.~\cite{SteinShakarchi}, extend
$g_u(x)=x h_t(x)e^{ux}$ for $x\ge0$ by $0$ for $x<0$.
The extension is continuous at $0$, is integrable, and decays
exponentially at $+\infty$. Its $\delta$-periodization
\[
 G_u(x)=\sum_{j\in\mathbb Z}g_u(x+j\delta)
\]
converges uniformly on $[0,\delta]$. Its $k$th Fourier
coefficient is
\begin{align*}
 \frac1\delta\int_0^\delta
  G_u(x)e^{-2\pi ikx/\delta}\,dx
 &=\frac1\delta\int_{\mathbb R}
    g_u(x)e^{-2\pi ikx/\delta}\,dx\\
 &=\frac1\delta H_t(u-2\pi ik/\delta).
\end{align*}
These coefficients are absolutely summable by the preceding
quadratic decay estimate. Fourier inversion at $x=0$, followed
by reindexing $k\mapsto-k$, gives
\begin{align}\label{eq:poisson-sampled}
 \delta\sum_{j=1}^\infty (j\delta)h_t(j\delta)e^{uj\delta}
 &=\sum_{k\in\mathbb Z}H_t(u+2\pi ik/\delta)\notag\\
 &=\delta^2e^{\delta u}F_{\delta,t}(e^{\delta u}).
\end{align}
Finally, when $|\Im u|\le\pi/\delta$ and $k\ne0$,
\[
 \left|\Im u+\frac{2\pi k}{\delta}\right|
 \ge\frac{(2|k|-1)\pi}{\delta}.
\]
Therefore
\begin{align*}
 \left|\sum_{k\ne0}H_t(u+2\pi ik/\delta)\right|
 &\le C_J\sum_{k\ne0}
   \frac1{1+\bigl((2|k|-1)\pi/\delta\bigr)^2}\\
 &\le\frac{2C_J\delta^2}{\pi^2}
    \sum_{k=1}^\infty\frac1{(2k-1)^2}
 \le C'_J\delta^2.
\end{align*}
Subtracting the $k=0$ term in \eqref{eq:poisson-sampled}
proves \eqref{eq:fixed-approx}.
\end{proof}

\begin{proof}[Proof of Theorem~\ref{thm:fixed}]
We separate the periodic part of the transform formula \eqref{eq:fixed-transform}:
\begin{equation}\label{eq:fixed-M-plus-error}
 u^2H_t(u)=M_t(u)
   -\frac{2(1-t)(e^u-1)}{u(1-te^u)},
 \qquad
 M_t(u)=\frac{B_t(e^u)}{(1-te^u)^2},
\end{equation}
where
\[
 B_t(y)=(1-ty)^2+(1-t)^2y
       =t^2y^2+(1-4t+t^2)y+1.
\]
Its discriminant and two endpoint values are
\begin{align*}
 \operatorname{disc}(B_t)&=(1-t)^2(1-6t+t^2)>0,\\
 B_t(-1)&=4t>0,\\
 B_t(-1/t)&=-\frac{1-6t+t^2}{t}<0.
\end{align*}
For $0<t<3-2\sqrt2$, the middle coefficient
$1-4t+t^2$ is positive, so the two roots of $B_t$
are negative. One of them is $-q$ with
\[
 1<q<t^{-1};
\]
it is the root of smaller modulus, since the product of the
moduli of both roots is $t^{-2}$. Consequently
$M_t$ is $2\pi i$-periodic and has the simple zeros
\[
 u_j=\log q+(2j+1)\pi i,\qquad j\in\mathbb Z,
\]
in $\Re u<L=-\log t$.

Choose
\[
 t<a<\min\{1,q^{-1}\},\qquad b=t/a,
\]
and then $r\in(0,\pi/2)$ sufficiently small that
\[
 0<\log q-r<\log q+r<-\log a.
\]
We may also require that the disk $\overline{D(u_0,r)}$
contains no zero of $M_t$ other than $u_0$. Periodicity yields
\[
 m_t:=\min_{|u-u_j|=r}|M_t(u)|>0,
 \qquad j\in\mathbb Z,
\]
with a minimum independent of $j$.
For $|u-u_j|=r$, $j\ge1$, we have $|u|\asymp j$.
Since $\Re u\in[\log q-r,\log q+r]\Subset(-\infty,L)$,
\eqref{eq:fixed-M-plus-error} gives
\[
 |u^2H_t(u)-M_t(u)|\le\frac{C_t}{j}.
\]
If the circle also lies in $|\Im u|\le\pi/\delta$,
Lemma~\ref{lem:fixed-periodization} gives
\begin{align}\label{eq:fixed-circle-error}
 &\left|u^2\delta^2e^{\delta u}
      F_{\delta,t}(e^{\delta u})-M_t(u)\right|\notag\\
 &\quad\le |u|^2
  \left|\delta^2e^{\delta u}F_{\delta,t}(e^{\delta u})-H_t(u)\right|
  +|u^2H_t(u)-M_t(u)|\notag\\
 &\quad\le C_t\bigl(\delta^2(j+1)^2+j^{-1}\bigr).
\end{align}
First choose $j_0$ with $C_t/j_0<m_t/3$.
Next choose $0<\eta<1/4$ with
$4C_t\eta^2<m_t/3$. For all sufficiently small $\delta$,
all disks with
\[
 j_0\le j\le\floor{\eta/\delta}
\]
lie in $0<\Im u<\pi/\delta$ and the right-hand side of
\eqref{eq:fixed-circle-error} is smaller than $m_t$.
Since $M_t$ has exactly one simple zero in each disk,
Lemma~\ref{lem:rouche} gives one simple zero of
$u^2\delta^2e^{\delta u}F_{\delta,t}(e^{\delta u})$
there. Each disk avoids $u=0$, and $\delta^2e^{\delta u}$
never vanishes, so this is exactly one simple zero of
$F_{\delta,t}(e^{\delta u})$.
The map $u\mapsto e^{\delta u}$ is injective in the strip
$0<\Im u<\pi/\delta$; it maps these disks into the upper
half-plane. The coefficients of \eqref{eq:fixed-samples} are
real for $t>0$, so conjugation supplies the same number of
distinct simple nonreal zeros below the real axis.

All the constructed zeros satisfy
\[
 0<\Re u<-\log a
 \quad\Longrightarrow\quad
 1<|e^{\delta u}|<a^{-\delta}.
\]
The product variables of $(av,a)$ and $(b,b)$ are $tv$ and $t$,
and
\[
 |av|^\gamma<a
 \quad\Longleftrightarrow\quad
 |v|<a^{1/\gamma-1}=a^{-\delta}.
\]
Both $a,b$ were chosen independently of $\delta$.
The resulting number of zeros is at least
\[
 2\bigl(\floor{\eta/\delta}-j_0+1\bigr)
 \ge\frac{\eta}{\delta}
 =\frac{\eta\gamma}{\gamma-1}
 \ge\frac{\eta}{\gamma-1}
\]
when $\delta$ is sufficiently small. This proves the lower
bound and the asserted simplicity.

For the upper bound, choose $\sigma$ with
$-\log a<\sigma<L$ and put $R=e^{\delta\sigma}$. Then
\[
 a^{-\delta}<R<t^{-\delta},
\]
so Jensen's formula can be applied on this circle, within the
section's domain. For $v=Re^{i\theta}$, $-\pi\le\theta\le\pi$, set
$u=\sigma+i\theta/\delta$. The Poisson approximation and the
decay of $H_t$ give
\begin{align*}
 \left|\delta^2e^{\delta u}F_{\delta,t}(e^{\delta u})\right|
 &\le |H_t(u)|+C_t\delta^2\\
 &\le\frac{C_t\delta^2}{\delta^2+\theta^2}+C_t\delta^2.
\end{align*}
As $|e^{\delta u}|=R$ stays bounded above and below,
\begin{equation}\label{eq:jensen-circle-bound}
 \left|F_{\delta,t}(Re^{i\theta})\right|
 \le C_t\left(1+\frac1{\delta^2+\theta^2}\right).
\end{equation}
The logarithmic mean on this circle is uniformly bounded:
\begin{align*}
 \frac1{2\pi}\int_{-\pi}^\pi
     \log^+|F_{\delta,t}(Re^{i\theta})|\,d\theta
 &\le C'_t+\frac1{2\pi}\int_{-\pi}^\pi
       \log\left(1+\frac1{\delta^2+\theta^2}\right)d\theta\\
 &\le C'_t+\frac1{2\pi}\int_{-\pi}^\pi
       \log\left(1+\frac1{\theta^2}\right)d\theta
 <\infty.
\end{align*}
Moreover
$F_{\delta,t}(0)=\beta_1+(1-\beta_1)t
 =1-(1-t)\delta\ge t$.
Every zero in $|v|<a^{-\delta}$ contributes at least
\[
 \log\frac{R}{a^{-\delta}}
 =\delta(\sigma+\log a)>0
\]
to Jensen's formula on $|v|=R$. Hence, writing $n_\gamma$
for their total multiplicity,
\begin{align*}
 n_\gamma\delta(\sigma+\log a)
 &\le\frac1{2\pi}\int_{-\pi}^\pi
       \log|F_{\delta,t}(Re^{i\theta})|\,d\theta
       -\log|F_{\delta,t}(0)|\\
 &\le C'_t-\log t.
\end{align*}
As $\delta^{-1}=\gamma/(\gamma-1)$, the upper bound
$n_\gamma\le C_t/(\gamma-1)$ follows upon decreasing
$\eta_t>0$ if necessary.
\end{proof}

\begin{proof}[Proof of Corollary~\ref{cor:fixed-pole}]
Set $t=1/10$ and recall the quadratic $B_t$ from the proof of
Theorem~\ref{thm:fixed}. Direct calculation gives
\begin{align*}
 B_{1/10}(-1)&=\frac25>0,\\
 B_{1/10}(-\sqrt{10})
 &=\frac{11}{10}-\frac{61\sqrt{10}}{100}<0.
\end{align*}
Consequently its smaller negative root is $-q$ with
$1<q<\sqrt{10}=t^{-1/2}$. We may therefore choose
\[
 a=b=\sqrt t=\frac1{\sqrt{10}},\qquad
 0<\log q<-\log a=\tfrac12\log10,
\]
in the proof of Theorem~\ref{thm:fixed}.
The fixed product variables are then
\[
 s=z_1\overline{w_{*,1}}=\frac{z_1}{\sqrt{10}},
 \qquad t=z_2\overline{w_{*,2}}=\frac1{10},
\]
which corresponds to $v=\sqrt{10}\,z_1$. The admissible disk
$|v|<a^{-\delta}$ is equivalent to
$|z_1|<10^{-1/(2\gamma)}$, giving the stated conclusion.
\end{proof}

\section{Proof of Theorem~\ref{thm:distribution}}\label{sec:distribution}

Recall that the degree of the limiting polynomial
$P_\delta$ is $d=N-1$. For a polynomial $G$ with $G(0)\ne0$ and a half-open
arc $I$ of the unit circle, let $n_G(I)$ count, with multiplicity,
the zeros $v$ of $G$ for which $v/|v|\in I$. We write $|I|$ for the
length of the arc in radians. We first estimate $n_{P_\delta}(I)$
uniformly in $I$ and $\delta$.

Hayman's sector estimate for sparse polynomials \cite{Hayman},
also stated in \cite[Theorem~1.3]{Hrubes}, asserts that a polynomial
$G$ of degree $D$ with $G(0)\ne0$ and at most $k$ nonzero terms satisfies
\begin{equation}\label{eq:hayman}
 \left|n_G(I)-D\frac{|I|}{2\pi}\right|\le k-1 .
\end{equation}
The half-open version follows by rotating both endpoint rays through a sufficiently small common angle and then taking a one-sided limit; there are only finitely many root arguments.

The coefficient of $v^j$ in $P_\delta(v)$ is
$(j+1)(1-\delta(j+1))$ for $0\le j\le d$. The generating
functions for $\sum_{j\ge0}(j+1)v^j$ and
$\sum_{j\ge0}(j+1)^2v^j$ give
\begin{align*}
 &\sum_{j=0}^\infty(j+1)(1-\delta(j+1))v^j\\
 &\qquad=\frac1{(1-v)^2}
      -\delta\frac{1+v}{(1-v)^3}
 =\frac{(1-\delta)-(1+\delta)v}{(1-v)^3}.
\end{align*}
Thus the third finite differences of the coefficients vanish in
the interior. Equivalently, the polynomial
$(1-v)^3P_\delta(v)$ has nonzero coefficients only at
$0,1,d+1,d+2,d+3$; this also follows from the explicit
identity \eqref{eq:sparse-limit}. It has degree $d+3$ and
at most five nonzero terms. Hayman's estimate therefore gives
\[
 \left|n_{(1-v)^3P_\delta}(I)
       -(d+3)\frac{|I|}{2\pi}\right|\le4.
\]
Since $(1-v)^3$ contributes exactly three zeros at $v=1$,
\begin{align*}
 n_{(1-v)^3P_\delta}(I)
   &=\begin{cases}
       n_{P_\delta}(I)+3,&1\in I,\\
       n_{P_\delta}(I),&1\notin I,
      \end{cases}\\
 \left|n_{P_\delta}(I)-d\frac{|I|}{2\pi}\right|
 &\le4+3\max\left\{\frac{|I|}{2\pi},
                            1-\frac{|I|}{2\pi}\right\}\le7.
\end{align*}
Hence
\begin{equation}\label{eq:angular-polynomial}
 \left|n_{P_\delta}(I)-d\frac{|I|}{2\pi}\right|\le7.
\end{equation}

\begin{lemma}\label{lem:radial}
Fix $0<\varepsilon<1$. For all sufficiently small $\delta>0$,
$P_\delta$ has no zeros in $|v|\le1-\varepsilon$ and at most one
zero, counted with multiplicity, in $|v|\ge1+\varepsilon$.
\end{lemma}

\begin{proof}
For each fixed $j\ge0$, the coefficient of $v^j$ in
$P_\delta$ tends to $j+1$ as $\delta\to0^+$, and
\[
 0\le(j+1)(1-\delta(j+1))\le j+1
 \qquad(0\le j\le d).
\]
For every $\rho<1$, the series $\sum_{j\ge0}(j+1)\rho^j$
converges. Dominated convergence for the power series thus gives
\[
 P_\delta(v)\longrightarrow
    \sum_{j=0}^\infty(j+1)v^j=\frac1{(1-v)^2}
 \quad\text{uniformly on }|v|\le\rho.
\]
The limit is zero-free. Taking $\rho=1-\varepsilon$ and
using Lemma~\ref{lem:rouche} (or Hurwitz's theorem) proves that
$P_\delta$ has no roots in $|v|\le1-\varepsilon$ for small
$\delta$.

For the exterior roots, write
$\delta^{-1}=N+\alpha$, $0<\alpha\le1$, and recall $d=N-1$.
Reversing the coefficients gives
\begin{align}\label{eq:reverse-distribution}
 v^dP_\delta(v^{-1})
 &=\sum_{j=0}^{d}
    (N-j)\left(1-\frac{N-j}{N+\alpha}\right)v^j\notag\\
 &=\sum_{j=0}^{d}\frac{(N-j)(\alpha+j)}{N+\alpha}v^j.
\end{align}
For $0\le j\le d$,
\[
 0\le\frac{(N-j)(\alpha+j)}{N+\alpha}\le \alpha+j\le j+1.
\]
Suppose that $\delta\to0^+$ along a sequence on which the
reversed polynomials have at least two roots in
$|v|\le(1+\varepsilon)^{-1}$. Passing to a subsequence, we may
assume $\alpha\to b\in[0,1]$. For each fixed $j$ the coefficient
in \eqref{eq:reverse-distribution} tends to $b+j$. The same
summable majorant gives locally uniform convergence on $|v|<1$:
\[
 v^dP_\delta(v^{-1})\longrightarrow
 \sum_{j=0}^\infty(b+j)v^j
 =\frac{b+(1-b)v}{(1-v)^2}.
\]
The limiting function has at most one zero in the unit disk.
Choose
\[
 (1+\varepsilon)^{-1}<\rho<1
\]
so that its numerator $b+(1-b)v$ does not vanish on $|v|=\rho$.
Lemma~\ref{lem:rouche} implies that the reversed polynomials have at
most one root in $|v|<\rho$ for all sufficiently small $\delta$,
a contradiction. Since the constant coefficient of $P_\delta$
is nonzero, inversion preserves root multiplicities and converts
$|v|\ge1+\varepsilon$ into
$|v^{-1}|\le(1+\varepsilon)^{-1}<\rho$. This proves the
second assertion.
\end{proof}

\begin{proof}[Proof of Theorem~\ref{thm:distribution}]
Fix $1<\gamma<2$ and let $d=\deg P_\delta$. Choose pairwise
disjoint disks around the distinct zeros of $P_\delta$, all
contained in $0<|v|<R_\gamma$, so small that each disk centered
at $\zeta$ satisfies
\begin{equation}\label{eq:root-disk-close}
 |v-\zeta|<\frac1d,\qquad
 \operatorname{dist}_{\mathbb R/2\pi\mathbb Z}
    (\arg v,\arg\zeta)<\frac1d.
\end{equation}
Every root is nonzero by $P_\delta(0)>0$, so these angle
conditions can be imposed simultaneously. By locally uniform
convergence $F_t\to P_\delta$, we can decrease $t_\gamma$
so that applying Lemma~\ref{lem:rouche} on each disk boundary gives the same
number of roots as its center's multiplicity, for all
$0<t<t_\gamma$. The large-disk count from
Theorem~\ref{thm:main} is still $d$; hence these disks contain
all roots counted there, and all of them are simple.

Let $I$ be an arc of length $|I|$. Enlarge $I$ by $1/d$ at each
endpoint, taking the whole circle if necessary. Similarly, shorten
$I$ by $1/d$ at each endpoint, taking the empty set if necessary.
By \eqref{eq:root-disk-close}, the root count of the section
inside $I$ lies between the polynomial counts for the shortened
and extended arcs. Their lengths differ from $|I|$ by at most
$2/d$. Applying \eqref{eq:angular-polynomial} to both arcs,
\[
 \left|\#\{\text{section roots with argument in }I\}
       -d\frac{|I|}{2\pi}\right|
 \le7+\frac{d}{2\pi}\frac2d
 =7+\frac1\pi<8.
\]
This estimate is uniform in the arc $I$. The interpretation for
half-open endpoints follows by one-sided approximation of the
endpoints, leaving all finitely many root arguments fixed.

Now fix $0<\varepsilon<1$. For sufficiently small $\delta$,
Lemma~\ref{lem:radial} with $\varepsilon/2$ shows that all but
at most one root of $P_\delta$ lie in
\[
 1-\frac\varepsilon2<|v|<1+\frac\varepsilon2.
\]
We also arrange $1/d<\varepsilon/2$. Then
\eqref{eq:root-disk-close} places all but at most one of the
matching section roots in
$1-\varepsilon<|v|<1+\varepsilon$.
Discarding that possible root changes the angular discrepancy by
at most one, so
\[
 \left|Z_{\gamma,t}(\varepsilon,I)
         -d\frac{|I|}{2\pi}\right|<9.
\]
Finally the definition of $d$ gives
\[
 d=\ceil{\frac1{\gamma-1}}-1,
 \qquad 0<\frac1{\gamma-1}-d\le1,
 \qquad 0\le\frac{|I|}{2\pi}\le1.
\]
Thus
\begin{align*}
 \left|Z_{\gamma,t}(\varepsilon,I)
       -\frac{|I|}{2\pi(\gamma-1)}\right|
 &\le\left|Z_{\gamma,t}(\varepsilon,I)
           -d\frac{|I|}{2\pi}\right|
     +\frac{|I|}{2\pi}\left|d-\frac1{\gamma-1}\right|\\
 &<9+1=10.
\end{align*}
All matching disks and the smallness condition on $t_\gamma$
are chosen before the arc $I$, as required.
\end{proof}


\par\medskip
\noindent\textbf{Data Availability Statement.}
No research datasets were generated or analyzed in this study.
\par\medskip
\noindent\textbf{Conflict of Interest.}
The author declares no competing interests.


\begin{thebibliography}{99}
\footnotesize
\linespread{1.0}\selectfont

\bibitem{Ahlfors1979}
L. V. Ahlfors,
\emph{Complex Analysis}, 3rd ed., McGraw-Hill, New York, 1979.

\bibitem{AhnPark2012}
H. Ahn and J.-D. Park,
The explicit forms and zeros of the Bergman kernel function for Hartogs type domains,
\emph{J. Funct. Anal.} \textbf{262} (2012), no.~8, 3518--3547.

\bibitem{Bell1982}
S. R. Bell,
The Bergman kernel function and proper holomorphic mappings,
\emph{Trans. Amer. Math. Soc.} \textbf{270} (1982), no.~2, 685--691.

\bibitem{BenderChakrabartiEdholmMainkar2022}
C. Bender, D. Chakrabarti, L. D. Edholm and M. Mainkar,
$L^p$-regularity of the Bergman projection on quotient domains,
\emph{Canad. J. Math.} \textbf{74} (2022), no.~3, 732--772.

\bibitem{Boas}
H. P. Boas,
Lu Qi-Keng's problem,
\emph{J. Korean Math. Soc.} \textbf{37} (2000), no.~2, 253--267.

\bibitem{BoasFuStraube1999}
H. P. Boas, S. Fu and E. J. Straube,
The Bergman kernel function: explicit formulas and zeroes,
\emph{Proc. Amer. Math. Soc.} \textbf{127} (1999), no.~3, 805--811.

\bibitem{ChakrabartiEdholm2024}
D. Chakrabarti and L. D. Edholm,
Projections onto $L^p$-Bergman spaces of Reinhardt domains,
\emph{Adv. Math.} \textbf{451} (2024), Art.~109790, 46 pp.

\bibitem{ChakrabartiEdholmMcNeal2019}
D. Chakrabarti, L. D. Edholm and J. D. McNeal,
Duality and approximation of Bergman spaces,
\emph{Adv. Math.} \textbf{341} (2019), 616--656.

\bibitem{ChakrabartiKonkelMainkarMiller2020}
D. Chakrabarti, A. Konkel, M. Mainkar and E. Miller,
Bergman kernels of elementary Reinhardt domains,
\emph{Pacific J. Math.} \textbf{306} (2020), no.~1, 67--93.

\bibitem{ChakrabartiZeytuncu2016}
D. Chakrabarti and Y. E. Zeytuncu,
$L^p$ mapping properties of the Bergman projection on the Hartogs triangle,
\emph{Proc. Amer. Math. Soc.} \textbf{144} (2016), no.~4, 1643--1653.

\bibitem{DengFu2026}
G. Deng and Q. Fu,
Bergman kernel zeros on rational Hartogs triangles,
preprint, arXiv:2611.07480 (2026).
\url{https://arxiv.org/abs/2611.07480}.

\bibitem{Edholm2016}
L. D. Edholm,
Bergman theory of certain generalized Hartogs triangles,
\emph{Pacific J. Math.} \textbf{284} (2016), no.~2, 327--342.

\bibitem{EdholmMathew2026}
L. D. Edholm and V. T. Mathew,
Arithmetic properties and zeros of the Bergman kernel on a class of quotient domains,
\emph{Math. Z.} \textbf{313} (2026), Art.~43, 27 pp.

\bibitem{EdholmMcNeal2017}
L. D. Edholm and J. D. McNeal,
Bergman subspaces and subkernels: degenerate $L^p$ mapping and zeroes,
\emph{J. Geom. Anal.} \textbf{27} (2017), no.~4, 2658--2683.

\bibitem{ErdosTuran1950}
P. Erd\H{o}s and P. Tur\'an,
On the distribution of roots of polynomials,
\emph{Ann. of Math.} (2) \textbf{51} (1950), no.~1, 105--119.

\bibitem{Hayman}
W. K. Hayman,
Angular value distribution of power series with gaps,
\emph{Proc. London Math. Soc.} (3) \textbf{24} (1972), no.~4, 590--624.

\bibitem{Hrubes}
P. Hrube\v{s},
On the real $\tau$-conjecture and the distribution of complex roots,
\emph{Theory Comput.} \textbf{9} (2013), no.~10, 403--411.

\bibitem{HuoWickWeak2020}
Z. Huo and B. D. Wick,
Weak-type estimates for the Bergman projection on the polydisc and the Hartogs triangle,
\emph{Bull. Lond. Math. Soc.} \textbf{52} (2020), no.~5, 891--906.

\bibitem{HuoWickWeighted2020}
Z. Huo and B. D. Wick,
Weighted estimates for the Bergman projection on the Hartogs triangle,
\emph{J. Funct. Anal.} \textbf{279} (2020), no.~9, Art.~108727, 34 pp.

\bibitem{SteinShakarchi}
E. M. Stein and R. Shakarchi,
\emph{Fourier Analysis: An Introduction},
Princeton Lectures in Analysis I, Princeton University Press,
Princeton, NJ, 2003.

\end{thebibliography}
\end{document}